\documentclass{article}

\usepackage[a4paper, left=3.3cm,top=3.3cm,right=3.3cm,bottom=3.7cm]{geometry}
\usepackage{concmath}
\usepackage[T1]{fontenc}

\usepackage{amsmath,amsfonts,amssymb}
\usepackage{graphicx}
\usepackage{mathtools}
\usepackage{framed}
\usepackage{enumitem}
\setenumerate{topsep=0em}
\setitemize{topsep=0em,label=\scriptsize{$\blacksquare$}}

\usepackage{siunitx}

\usepackage{authblk}
\usepackage{tikz}
\usepackage{tikz-3dplot}

\newcommand{\edot}{\,\cdot\,}

\newcommand*{\R}{\mathbb{R}}

\newcommand\abs[1]{\left\vert#1\right\vert}

\newcommand\norm[1]{\left\Vert#1\right\Vert}
\newcommand\snorm[1]{\Vert#1\Vert}

\newcommand\set[1]{\bigl\{#1\bigr\}}
\newcommand{\kl}[1]{\left(#1\right)}

\newcommand\inner[2]{\left\langle#1,#2\right\rangle}

\DeclareMathOperator{\argmin}{\mathbf{argmin}}
\DeclareMathOperator{\prox}{\mathbf{prox}}
\DeclareMathOperator{\Ro}{\mathbf{R}}
\DeclareMathOperator{\Xo}{\mathbf{X}}
\DeclareMathOperator{\Wo}{\mathbf{W}}

\DeclareMathOperator{\ft}{\mathcal{F}}

\DeclareMathOperator{\Lambdaop}{\boldsymbol{\Lambda}}

\newcommand{\rr}{\mathbf{r}}

\newcommand{\source}{h}
\newcommand{\psource}{H}

\allowdisplaybreaks

\title{Photoacoustic image reconstruction from full field data in heterogeneous media}

\author{Markus Haltmeier}

\affil{Department of Mathematics, University of Innsbruck\authorcr
Technikerstrasse 13, 6020 Innsbruck, Austria\authorcr
E-mail: {\tt markus.haltmeier@uibk.ac.at}}

\author{Gerhard Zangerl}

\affil{Department of Mathematics, University of Innsbruck\authorcr
Technikerstrasse 13, 6020 Innsbruck, Austria\authorcr
E-mail: {\tt gerhard.zangerl@uibk.ac.at}}

\author{Robert Nuster}

\affil{Department of Physics, Universit\"at Graz\authorcr
Universitaetsplatz 5, Graz, Austria\authorcr
E-mail: {\tt ro.nuster@uni-graz.at }}

 \author{Linh V. Nguyen}

\affil{Department of Mathematics, University of Idaho\authorcr
875 Perimeter Dr, Moscow, ID 83844, US\authorcr
E-mail: {\tt lnguyen@uidaho.edu}}

\date{Januar 19, 2019}

\begin{document}
\maketitle

\begin{abstract}
We consider image reconstruction in  full-field photoacoustic tomography, where 2D projections  of the full  3D acoustic pressure distribution
at a given time $T>0$  are  collected.
We discuss existing results on the stability and uniqueness of
the resulting image reconstruction  problem and review existing reconstruction algorithms. Open challenges  are also mentioned.
Additionally,  we introduce   novel one-step  reconstruction methods
allowing for a variable speed of sound.  We apply preconditioned  iterative and variational regularization methods to the one-step formulation. Numerical results using the one-step formulation are presented,  together with a comparison  with   the previous  two-step approach for
full-field photoacoustic tomography.

 \end{abstract}

\section{Introduction}
\label{sec:intro}

Photoacoustic tomography (PAT) is a hybrid imaging modality that beneficially combines the high spatial resolution of ultrasound imaging with the  good  contrast of optical tomography \cite{wang2011photoacoustic,Wan09b}.
In PAT, a semi-transparent sample is illuminated by a short laser pulse which induces  an acoustic pressure wave depending
on the light absorbing structures inside  the sample.
The induced pressure waves propagate in space,  are
 detected outside of the  sample, and measurements are used
 to recover the photoacoustic (PA) source.
The standard approach in PAT is to record time-resolved
acoustic signals  on a detection surface  partially or fully enclosing the
investigated object. Plenty of reconstruction methods
have been
developed for the  standard setting. This
includes analytic inversion methods \cite{IPI,finch2007inversion,FPR,Halt-Inv,Kun07,XW05,Burgholzer2007},
time reversal  \cite{US,qian2011new,HKN}, continuous
iterative methods \cite{arridge2016adjoint,belhachmi2016direct,haltmeier2017analysis} and  discrete iterative methods \cite{wang2012investigation,DeaBueNtzRaz12,huang2013full}.

In this paper, we consider full field detection PAT (FFD-PAT),
which collects  different measurement data. In  FFD-PAT,
2D  linear projections   of the 3D pressure field   at a fixed  time instant $T>0$ are captured \cite{nuster2010full}.   This can be  implemented
by using a  special phase contrast  method  and a CCD-camera that records the full field projections \cite{nuster2014high}.
Similar to the  integrating line detector approach. \cite{burgholzer2007temporal,GruEtAl10,PalNusHalBur07a}, a 3D data  set is   measured    by collecting 2D full field projections
from a 1D set  of projection directions. Previous reconstructions
techniques for FFD-PAT where mostly based on a constant sound speed  assumption \cite{nuster2010full,nuster2014high}.
Recently, in \cite{zangerl2018full} we
investigated the case of a variable  speed of sound and
 established  the  stability  and  uniqueness
 in  a complete data situation. Moreover, we  proposed a  two-step reconstruction  procedure (including the partial data  case), where in the first step the 3D  pressure field
$p(\edot,  T)$ is recovered using  the measured projection data.
In the second step,    the  PA source     $\source = p(\edot,  0)$   is  recovered from $p(\edot,  T)$ by solving a finite time wave inversion problem.

We study image reconstruction in FFD-PAT allowing for a spatially variable speed of sound. In section \ref{sec:full} we   formulate the reconstruction  problem and recall  existing results.   In Section \ref{sec:one}, we present the new one-step formulation  where we recover the PA source directly from  the full field data. To stabilize  the  inversion, we use variational regularization (generalized Tikhonov  regularization)
including a preconditioning strategy   and minimize the  Tikhonov functional by  proximal  forward-backward
splitting. Numerical results for the one-step and the two-step method are
 presented  in Section~\ref{sec:num}.
 The paper ends with  some conclusions  given in   Section~\ref{sec:conclusion}.

\section{Full field detection in photoacoustic tomography}
 \label{sec:full}

In this section, we  formulate the
image reconstruction  problem in FFD-PAT and review known  uniqueness and stability results as well as
existing image reconstruction methods.

\subsection{Modeling and problem formulation}
 \label{sec:model}

We allow  a variable sound speed and model the
acoustic wave propagation in PAT by the
wave equation
\begin{align}\label{eq:wave1}
\kl{ \partial_t^2    - c^2(\rr) \Delta } p(\rr, t) &= 0,
&&
\text{ for } (\rr, t) \in \R^3 \times (0, \infty) \\ \label{eq:wave2}
\kl{ p(\rr, 0 ), \partial_t p(\rr, 0) } &= \kl{ \source(\rr),0}
&&
\text{ for }  \rr  \in \R^3  \,.
\end{align}
Here $c(\rr) > 0$ is the sound speed at location $\rr \in \R^3$ and $\source \colon \R^3 \to \R$ is the  initial
pressure distribution (the PA source) that encodes the inner
structure of the sample. We  assume that the PA source vanishes outside a bounded volume  $V  \subseteq \R^3$  and that the sound speed takes  the constant value $c_0$  outside a bounded volume containing $V$.

To formulate  the FFD-PAT reconstruction problem we
introduce  some   notations. We denote by   $p(\edot,T) = \Wo_T \source  $ the solution  of the wave equation   \eqref{eq:wave1}, \eqref{eq:wave2}  with initial data $\source$ at  the given measurement
time $T>0$, and  the initial-to-final time wave operator by
\begin{equation}\label{eq:Wo}
\Wo_T  \colon C_0^{\infty} ( V )
\to C_0^{\infty} (V_T) \colon
\source  \mapsto p(\edot,T) \,.
\end{equation}
Here and below  $C_0^{\infty} ( V )$ is the set of all smooth functions  in $\R^3$ that have compact support in $V$; the same  notation  is used for
$C_0^{\infty} (V_T)$, where  $V_T$ is sufficiently large. We define the  X-ray transform
\begin{equation}\label{eq:Radon}
\Xo \colon   C_0^\infty \left(  V_T  \right)  \rightarrow L^2 ( S^1   \times \R^2  )\,,
\quad (\Xo \phi)( \theta,  \xi, z) \coloneqq
\int_{\R} \phi(       \xi  \theta + s \theta^\bot , z   ) \, ds \,.
\end{equation}
Then $\Xo \source (\edot,  \edot, z) $ equals the  2D
Radon transform of $\source$ applied in horizontal
planes  $z= \text{const}$.    In particular,  $\Xo$ can be inverted using
 any of the existing reconstruction  formulas for  the
 2D Radon transform \cite{Dea83,Hel99,natterer2001mathematics}.

In the  complete data situation,  the challenge in FFD-PAT
 is to reconstruct the  PA source  $\source$ from
projection data $ P_T(\theta, \edot)  = \Xo  \Wo_T  \source$.
In practice,  one cannot measure  the function
$ P_T(\theta, \edot) $  on the whole projection plane.
In particular,  at least  integrals over  lines
intersecting  the investigated object are missing.
 Due  to practical constraints,  projection data
 might be unknown  for even more lines \cite{nuster2010full,nuster2014high}.
 Thus, in FFD-PT  we face with  the   following
 reconstruction problem.

\noindent\textbf{Problem~1 (Image reconstruction problem in FFD-PAT).}
{\em Let   $T>0$ be a given final measurement  time and
for any  $\theta \in S^1$, let  $ M_\theta \subseteq \R^2$
be the  set in the projection plane  where  the projection data are available.
The goal in FFD-PAT is to recover the PA source
 $\source$ from  data $P_{T,M}  = \chi_M  ( \Xo  \Wo_T  \source)$, where
\begin{equation}\label{eq:fm}
\chi_M   \Xo  \Wo_T  \colon C_0^{\infty} ( V ) \to
L^2(   S^1   \times  \R^2 ) \colon  \source
\mapsto \chi_M  ( \Xo  \Wo_T  \source)  \,,
\end{equation}
(with $\chi_M (\theta , \xi, z) =  1$ for $\xi \in M_\theta$ and
$\chi_M (\theta , \xi, z) =  0$ for $\xi \in \R^2 \setminus M_\theta$)
is the FFD-PAT forward operator.}

Problem~1 is closely related to the following finite time wave inversion
problem, where in particular the spatial dimensions $d = 1,2,3$ appear in
FFD-PAT.

\noindent\textbf{Problem~2 (Finite time wave inversion  problem).}
{\em Let $d \geq 1$ denote the spatial dimension.
For given   $T>0$  and a region $ D \subseteq \R^d$ where
pressure measurements are made,       recover the source  $\source$ from  data
$ \chi_D  \Wo_{T}  \source$
where $\Wo_{T}  \source$ is the solution
of the $d$-dimensional wave equation $\kl{ \partial_t^2    - c^2(\rr)  \Delta } p =0$ with  initial conditions $( p, \partial_t p )|_{t=0} = (\source,0) $.
}

In the remainder of this section we review known
results and open problems related  to Problems 1 and 2.
Especially in the case of a variable sound speed,
not so much is  known for these image reconstruction problems.

\subsection{Review  for  the constant sound speed case}

Problem~1 has  first been studied in \cite{nuster2010full} for    a  constant speed of sound $c=c_0$, where we  proposed  a non-iterative reconstruction method that is outlined in the following.
The commutation relation between  the 2D Radon transform and the wave equation states   $\chi_M \kl{ \Xo  \Wo_T \source } (\theta, \edot)=  \chi_M \kl{  \Wo_T    \Xo \source (\theta, \edot) } $.  This shows that
the FFD data, for any $\theta$, are given by  $\chi_{M_\theta}
P_ \theta(\edot, T )$ where
$P_ \theta$ solves  the 2D wave  equation
\begin{align}\label{eq:wave2d-1}
\kl{ \partial_t^2    - c_0^2 \Delta } P_\theta(\xi, z, t) &= 0,
&&
(\xi, z, t) \in \R^2 \times (0, \infty)
\\ \label{eq:wave2d-2}
\kl{ P_\theta(\xi, z, 0 ), \partial_t P_\theta(\xi, z, 0)  }
&=  \kl{(\Xo \source)(\theta, \xi, z), 0},
&&
(\xi, z) \in \R^2   \,.
\end{align}
Hence the FFD-PAT Problem~1 amounts to the solution
of the  2D instance of Problem~2 for every $\theta$.
If we have uniqueness, stability or an inversion
procedure for  Problem~2, we also have corresponding
results for Problem~1.

In the general case, uniqueness and stability for  Problem  2
are unknown.  For some special cases, in \cite{nuster2010full}
we derived an inversion  method  for the 2D instance based  on the  reduction to the  1D instance of Problem~2.
Below we formulate this  procedure for arbitrary dimension.
For that purpose, consider  the case of full data $M = \R^d$
and denote by
\begin{equation} \label{eq:radon}
\Ro \colon
C_0^\infty(V_T) \to L^2  (S^{d-1} \times \R)
\,, \quad
\Ro \phi ( \alpha , u  ) \triangleq
\int_{ \alpha^\bot } \phi( u \alpha + w  ) d w \,,
\end{equation}
the $d$-dimensional Radon transform that maps a  function
$\phi $ defined on $\R^d$ and vanishing outside $V_T$ to the integrals
of $f$ over all hyperplanes ($(d-1)$-dimensional affine subspaces) of $\R^d$.
Let $p $ denote the solution of the $d$-dimensional wave equation
with   initial conditions $( p, \partial_t p )|_{t=0} = (\source,0) $.
The  commutation relation between  the Radon transform
 and the wave equation  shows that
  $Q_\alpha  \triangleq   \Ro p (\alpha, \edot)$
 satisfies the 1D wave  equation $\kl{ \partial_t^2    - c^2 \partial_u^2 } Q_\theta = 0$ with initial conditions
$(Q_\alpha, \partial_t Q_\alpha)(\edot , 0 )   =(\Ro \source (\alpha, \edot), 0) $.
Evaluating the solution of the 1D  wave equation
at   $T$ gives   $ \Ro p (\alpha, s , T ) =  \frac{1}{2}
(
\Ro  \source ( \alpha, s - c T  )
+
\Ro  \source ( \alpha, s + c T  )
) $. If $ \source$ vanishes outside the ball $B_a =
\set{\rr\mid \abs{\rr} < a}$, then
 $\Ro  \source  ( \alpha, s ) = 0$ for $\abs{s} \geq a$.
  For $ cT \geq   a$  this implies
 \begin{equation}
 \Ro  p ( \alpha, s, T ) =
 \frac{1}{2}
 \begin{cases}
\Ro  \psource ( \alpha, s - c T  )  & \text{ for } s \geq  0 \,,
\\
\Ro  \psource ( \alpha, \xi + c T  )  & \text{ for } s \leq  0 \,.
  \end{cases}
\end{equation}
In particular,   $ \Ro  p ( \alpha, \edot , T ) $ consists of two  separated and translated  copies  of $\Ro  \source  ( \alpha, \edot)$, where each of the copies can be used to recover the 1D initial source $\Ro  \source
 ( \alpha, \edot)$.
This implies the following  results.

\noindent\textbf{Theorem 3.}
{\em If $ cT \geq   a$, $\source(\rr)$ vanishes
outside $B_a$ and $ D \supseteq  B_{a+cT}$,
then  the final time   wave inversion
Problem~2 is uniquely solvable.
Moreover, the initial source  $\source$ can be reconstructed from
$g  = \Wo_{T} \source$ by  the following procedure
\begin{enumerate}[label=(W\arabic*), itemsep=0em, leftmargin=3em]
\item Compute the Radon transform $G = \Ro g$.
\item Invert the 1D wave equation:
Choose $w_{-1} + w_1=1$, set $H (\alpha, \xi) \triangleq  \sum_{\sigma=\pm1} w_\sigma G ( \alpha, \xi +\sigma cT)$  for   $\abs{\xi} \leq a$.

\item
Compute   $\source  =  \Ro^{-1}  H$  where $\Ro^{-1} $
is any reconstruction  method  for the  Radon transform.
\end{enumerate}
}

Uniqueness in Theorem 3 refers to the fact that $\Wo_{T} h_1
\neq \Wo_{T} h_2$ whenever $h_1 \neq h_2$ are
two distinct PA sources.     Moreover, using  the two-sided stability  estimate
 $\norm{ \source}_{L^2} \asymp   \norm{ \Ro \source}_{H^{(d-1)/2}}$
 for the Radon transform,  the  reconstruction approach
 of Theorem 3 implies the two-sided stability estimate
 with respect to the $L^2$-norms
 \begin{equation} \label{eq:stabW}
\forall \source \in C_0^\infty(B_a) \colon \quad
\norm{ \source}_{L^2}
 \asymp  \norm{ \Wo_T \source}_{L^2}      \,.
 \end{equation}
Here and elsewhere, the notion $\norm{ \source}_{L^2}
 \asymp  \norm{ \Wo_T \source}_{L^2}  $ means that the inequalities
 $ c_1 \norm{ \source}_{L^2}
\leq
\norm{ \Wo_T \source}_{L^2}
\leq
c_2 \norm{\source}_{L^2} $ hold
for  some constants  $c_1, c_2 < \infty$.
Combining this with the considerations at the beginning
of this subsection we obtain the following result for  FFD-PAT.

\noindent\textbf{Theorem 4.}
{\em If $V = B_a$, $M_\theta \supseteq  B_{a+cT}$ and $ cT \geq   a$,
then the FFD-PAT Problem~1 is uniquely
solvable via $ \source = \Xo^{-1}  \Wo_T^{-1} G$.
Moreover the   two-sided  stability estimate
$\norm{ \source }_{L^2}  \asymp \norm{  \Xo  \Wo_T   \source  }_{H^{1/2}}$ holds. }

In \cite{nuster2010full}  we also apply the  Radon
transform approach to the limited data
case  $M_\theta \supsetneq  B_{a+cT}$.

\subsection{Review  for  the variable sound speed case}

In the case of  variable speed of sound where $c(\rr)$ is not
constant, the X-ray projections $\Xo  p (\theta, \edot)$ do
 not satisfy the 2D wave  equation. In \cite{zangerl2018full}
we therefore proposed  a  different approach where
we first invert the X-ray transform  and subsequently
solve the 3D finite time wave inversion Problem~2.

The following non-trivial fact  for Problem~2 has been proven in  \cite{zangerl2018full}.

\noindent\textbf{Theorem 5.}
{\em  For every $T>0$, the  finite time wave inversion Problem~2
with $D = \R^d$ is uniquely solvable. Moreover, the two-sided stability estimate  $\norm{ \source }_{L^2}  \asymp \norm{    \Wo_T   \source  }_{L^2}$ holds.}

As a consequence we have the following  result for
Problem~1 in the variable sound speed case.

\noindent\textbf{Theorem 6.}
{\em  For every $T>0$, the FFD-PAT
Problem~1 with $M_\theta  = \R^2$  is uniquely
solvable via $ \source =  \Wo_T^{-1} \Xo^{-1} G$.
Moreover, the two-sided stability estimate  $\norm{ \source}_{L^2}
\asymp \norm{ \Xo  \Wo_T   \source }_{H^{1/2}}$ holds.
}

For  the limited data case, in \cite{zangerl2018full}   we
proposed  the two-step procedure
 $\Wo_T^{-1} \Xo^{-1} (\chi_M \Xo  \Wo_T   \source )$.
In the case that $T$  is sufficiently
large we observed accurate reconstruction results using
the two-step  procedure.   However, $\Xo^{-1}$ does not exactly invert
$\chi_M \Xo$ and therefore the two-step method introduces a systematic (albeit small) error. Moreover, let us mention  that in the limited-data case, no  theoretical results on the  uniqueness and  stability  for Problem~1  are known.

\section{Preconditioned one-step inversion methods}
\label{sec:one}

The application of the inverse X-ray  transform $\Xo^{-1}$
in the two-step method for  limited data introduces  a systematic
error. Therefore, in this paper we propose a
one-step method  where  we recover $\source$ directly from
data $\chi_M \Xo  \Wo_T   \source$ instead of first
applying the  inverse  X-ray transform.
For implementing  the one-step strategy we
use iterative and variational  reconstruction methods.
We include a  preconditioned technique accounting for the
smoothing of  $\Xo$ by degree $1/2$ in order to  accelerate
the iteration.

\subsection{Preconditioned one-step   Landweber  method}

Denote by $\Lambdaop $
the preconditioning operator defined by 
 $\ft_2 \Lambdaop  \Phi (\theta, \omega, z)   \triangleq
 \frac{\abs{\omega}}{4\pi}  (\ft_2 \Phi)
  (\theta, \omega, z) $, and define  the backprojection  $\Xo^* \Phi (x,y,z)
\triangleq
\int_{S^1}  \Phi(\theta, \inner{\theta}{(x,y)}, z) \, d\alpha $
for $\Phi \colon S^1 \times \R^2 \to \R$.
Here  $\ft_2$ is the  Fourier transform in second  variable.
The backprojection operator  is the formal $L^2$-adjoint
of $\Xo$ with the respect to the standard $L^2$-inner product.
The  composition $\Xo^*  \Lambdaop= \Xo^{-1} $ is the standard  filtered   backprojection inversion formula for the 2D Radon transform  \cite{natterer2001mathematics} applied with fixed $z$.
We write $ \Wo_T^* $ for the formal
  $L^2$-adjoint of  the initial-to-final time  wave operator.
In \cite{zangerl2018full}  we have shown that
$\Wo_T^*  g =  \chi_{V}(\edot) \, q(\edot, 0)$ where   $q$
is the solution of the backwards wave equation
$( \partial_t^2 q  - c(\rr)^2  \Delta ) q  = 0 $ on
 $ \R^3 \times (-\infty,T)$  with
  $ (q , \partial_t q)|_{t=T}       = ( g, 0)$
  and $\chi_{V} $ denoting the indicator function of  $V$.

\noindent\textbf{Algorithm 7 (Preconditioned one-step     Landweber algorithm for  FFD-PAT).}{\em
\mbox{}
	\begin{enumerate}[label=(S\arabic*), itemsep=0em]
		\item  Initialize: $k = 0$,  $\source_0 = 0$,
		\item  While (not stop) do:
		\begin{itemize}
		  \item $ G_k   =  \Xo \Wo_T \source_k $.
		  \item $ r_k   =  \Wo_T^* \Xo^* \Lambdaop(G_k-G)    $.
		  \item $ \source_{k+1} =  \source_{k} - s_k  r_k  $
		  for step size $s_k  >0$
		\end{itemize}
		\end{enumerate}	
}

 Note that $\Wo_T^* \Xo^* \Lambdaop $ is the adjoint of the forward operator with respect  the weighted inner product $\inner{\Lambdaop  G_1}{G_2}_{L^2}$. Hence the above algorithm  undoes   the implicit smoothing of   $\Xo$ and $\Xo^*$. In particular, assuming the
 stability estimate  $\norm{ \source }_{L^2}  \asymp \norm{    \Wo_T   \source  }_{L^2}$, Algorithm 7 is linearly convergent:
 \begin{equation}
 \norm{  \source_{k+1} - \source_\star }_{L^2} \leq
  c \norm{  \source_k - \source_\star }_{L^2}
  \quad \text{ for some constant } c < 1 \,,
 \end{equation}
 where $G = \chi_M \Xo  \Wo_T   \source_\star$. While we expect the stability estimate  $\norm{ \source }_{L^2}  \asymp \norm{    \Wo_T   \source  }_{L^2}$ to be satisfied for important  special cases,
  it is expected not to be satisfied for  severe limited data cases  or a trapping sound  speed. In this case, we  expect Problems 1 and 2  to be  severely ill-posed, and that  no estimate of the form
  $\norm{ \source }_{L^2}
 \leq c \norm{    \Wo_T   \source  }_{H^\alpha} $ holds.
 Variational methods that include  an additional  regularization term
 are a  reasonable alternative in  the ill-posed and the
 well-posed case.

\subsection{Preconditioned one-step variational regularization}

Instead of looking  for a theoretically exact solution of the one-step formulation,
in  variational regularization  we minimize   the  generalized
Tikhonov functional
\begin{equation}
\mathcal{T}_{g, \lambda} (\source) \triangleq
 \frac{1}{2} \norm{ \Lambdaop (\chi_M \Xo  \Wo_T \source - G) }_{L^2}^2 + \lambda   \mathcal{R} (\source) \,,
 \end{equation}
 where $\mathcal{R}$ is a regularization term, that additionally stabilizes
 the iteration and $\lambda   \geq 0$ is the regularization parameter.
For its  numerical solution, we use the forward backward splitting \cite{ComWaj05}, which alternates between a  gradient step for the data fitting   term
$ \frac{1}{2} \norm{ \Lambdaop (\chi_M \Xo  \Wo_T \source - G) }_{L^2}^2$ and an  implicit step with respect to the regularization term $\lambda  \mathcal{R} (\source)$.

\noindent\textbf{Algorithm 8 (Preconditioned proximal
one-step algorithm  of FFD-PAT).}{\em
\mbox{}
	\begin{enumerate}[label=(S\arabic*), itemsep=0em]
		\item  Initialize: $k = 0$,  $\source_0 = 0$,
		\item  While (not stop) do
		\begin{itemize}
		  \item
		  $ r_k   =  \Wo_T^* \Xo^* \Lambdaop(\Xo \Wo_T \source_k-G)    $.
		  \item $ \source_{k+1} =  \prox_{s_k \lambda  \mathcal{R} }
		  \kl{ \source_{k} - s_k  r_k } $ for step size $s_k > 0$.
		\end{itemize}
		\end{enumerate}	
}

\bigskip
The  proximal mapping $\prox_{ s_k \lambda  \mathcal{R} }$ in
  Algorithm 8  is defined by
\begin{equation}
\prox_{ s_k \lambda  \mathcal{R} } (   f ) =  \argmin_{\source \in L^2 (V)} \left\{ \frac{1}{2} \norm{ \source - f }_{L^2}^2 +
s_k \lambda  \mathcal{R}  (\source)   \right\}\,.
\end{equation}
This  implicit treatment of the regularizer $\mathcal{R}$ allows efficient treatment of  non-differentiable regularizers,
where  the plain gradient method is not applicable.
Also in the   differentiable case,  the forward-backward splitting   is
useful. For example, if  $\mathcal{R}(\source)  = \frac{1}{2} \snorm{\nabla \source}_{L_2}^2$, then the Tikhonov functional  $\mathcal{T}_{g, \lambda}$  is ill-conditioned. The   plain  gradient  method therefore becomes inefficient  whereas the
  forward-backward splitting treats  $\frac{1}{2} \snorm{\nabla \source}_{L_2}^2$ implicitly and is not affected by the ill-conditioning.
  An additional benefit of  Algorithm 8   is the  explicit smoothing step, which we observed to overall stabilize the  iterative procedure.

 \section{Numerical results}
\label{sec:num}

For the presented numerical results, we  consider  a 2D version
of the FFD-PAT Problem~1. This arises when we have translational symmetry in the $z$-direction meaning that  the sound speed
and the PA source  are independent of the $z$-direction.
We assume   that the PA source $\source$ is  contained
in the 2D  ball $D_a$. The 2D full field data
are given by   $ G = \chi_{M}  \Xo \Wo_T \source$
where  $\Wo_T$ is the solution of the 2D wave   equation at time $T$ and $\Xo$ reduces to the  2D Radon transform.
The   measurement  region is taken as
$M_\theta =  \R \setminus [-a,a]$  for any $\theta \in S^1$.

\subsection{Implementation details}

All  numerical  results are obtained using  Matlab.
In  the numerical implementation we  solve the  2D wave equation
and its adjoint using the  k-space method \cite{mast2001k,cox2007k,tabei2002k} as described
in \cite{haltmeier2017analysis}.
The Radon   transform, its adjoint and inverse we computed with   the Matlab build in functions \texttt{radon} and  \texttt{iradon} using linear interpolation.
 Data simulation consist in first  computing the solution
 of the wave equation
 $p = \Wo_T \source $, then computing  the X-ray transform
 $P = \Xo p$ and finally multiplying with   $\chi_ {M_\theta}$ where
 $M_\theta =\R \setminus [-a,a]$ with $a=1$.
 We  discretize the PA source  using  a $201 \times 201$ Cartesian
  grid on the domain $[-1,1]^2$ and numerically compute
  $ \Wo_T \source $ at time $T=2$ on an $601 \times 601$ Cartesian
  grid on the domain $[-3,3]^2$.
  The X-ray transform is evaluated  for $1000$ directions  equidistantly distributed between $0$ and $\SI{180}{\degree}$.

For image reconstruction   we use the proximal gradient version of the
one-step method (Algorithm 8). For comparison purpose we also apply  the  two-step method   $p = \Wo_T^{-1} \Xo^{-1} G$.
For stably inverting $\Wo_T$ we use the  proximal gradient
method for minimizing the      Tikhonov functional $\frac{1}{2} \norm{   \Wo_T \source - G }_{L^2}^2 + \lambda   \mathcal{R} (\source)$. This results in the following two-step variant of Algorithm~8.

\begin{figure}[htb!] \centering
\includegraphics[width=0.48\textwidth]{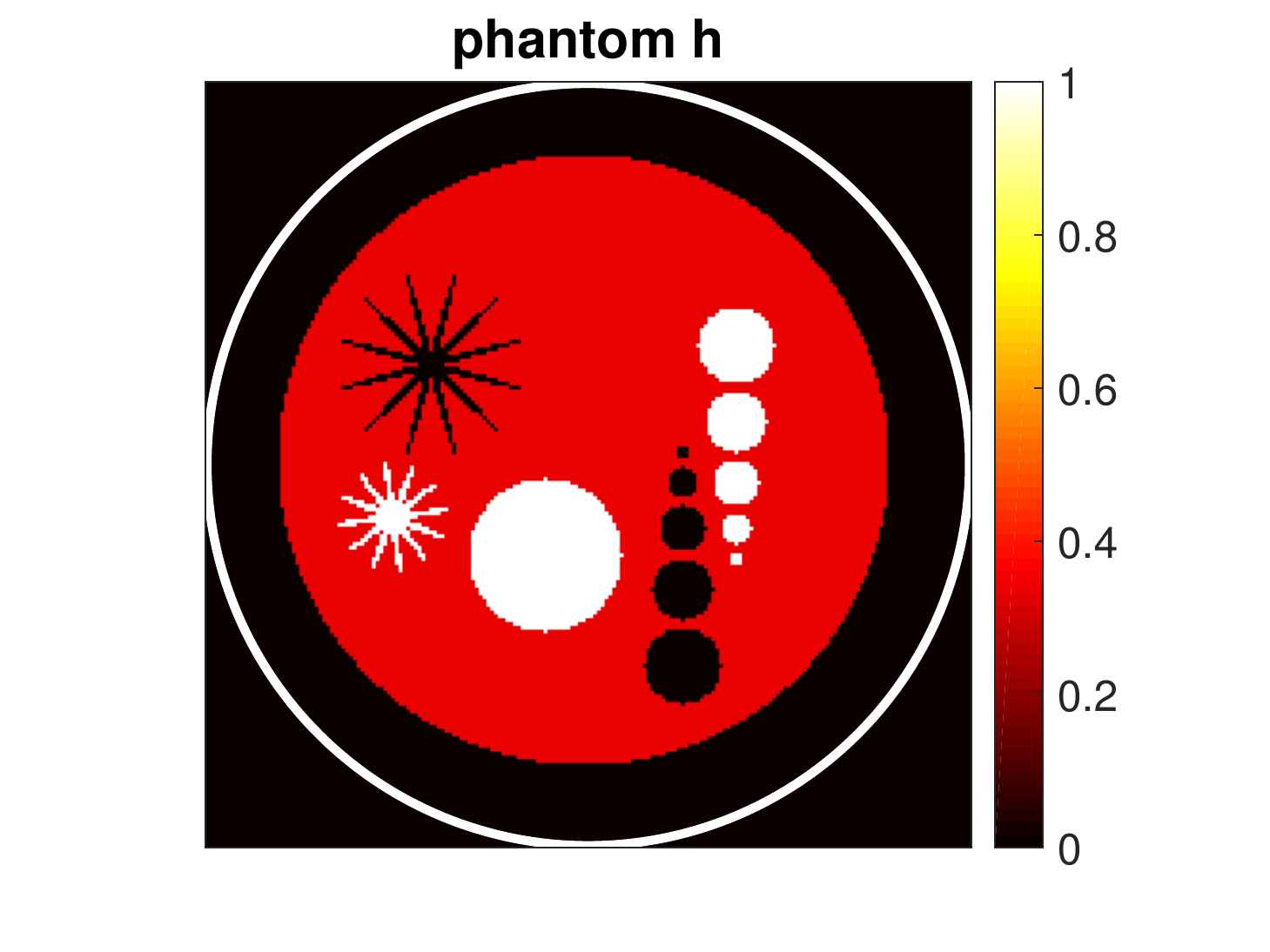} \;
\includegraphics[width=0.48\textwidth]{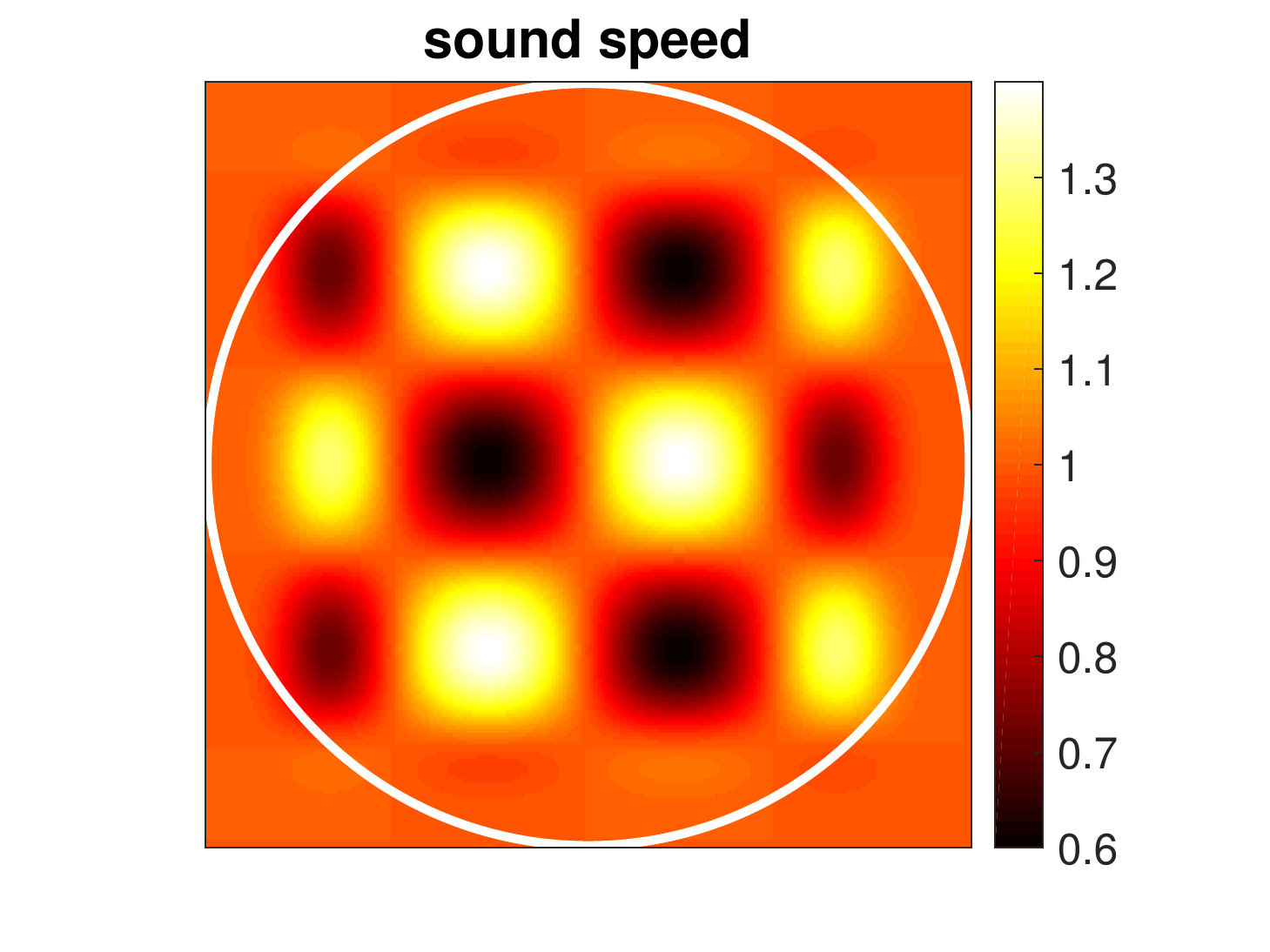} \;
\caption{ \textbf{Phantom and sound speed.}  Left: PA source
 to be reconstructed. Right: Used  trapping sound speed
 profile taken from  \cite{qian2011new}.}
\label{fig:phantom}
\end{figure}

\noindent\textbf{Algorithm 9 (Proximal
two-step algorithm  for FFD-PAT).}{\em
\mbox{}
	\begin{enumerate}[label=(S\arabic*), itemsep=0em]
		\item Compute $g = \Xo^{-1} G    $.
		\item  Initialize: $k = 0$,  $\source_0 = 0$,
		\item  While (not stop) do:
		 $ \source_{k+1} =  \prox_{s_k \lambda  \mathcal{R} }
		  \kl{ \source_{k} - s_k   \Wo_T^*  (\Wo_T \source_k  - g) } $
		  for step size $s_k > 0$.
		\end{enumerate}	
}

For the one-step and the two-step  algorithm
we use the regularizer $\mathcal{R} (\source) = \frac{1}{2} \snorm{\nabla \source}_{L_2}^2$, the regularization parameter
$\lambda  =0.5$, and  a constant step size $s_k =0.2$.

\begin{figure}[htb!] \centering
\includegraphics[width=0.48\textwidth]{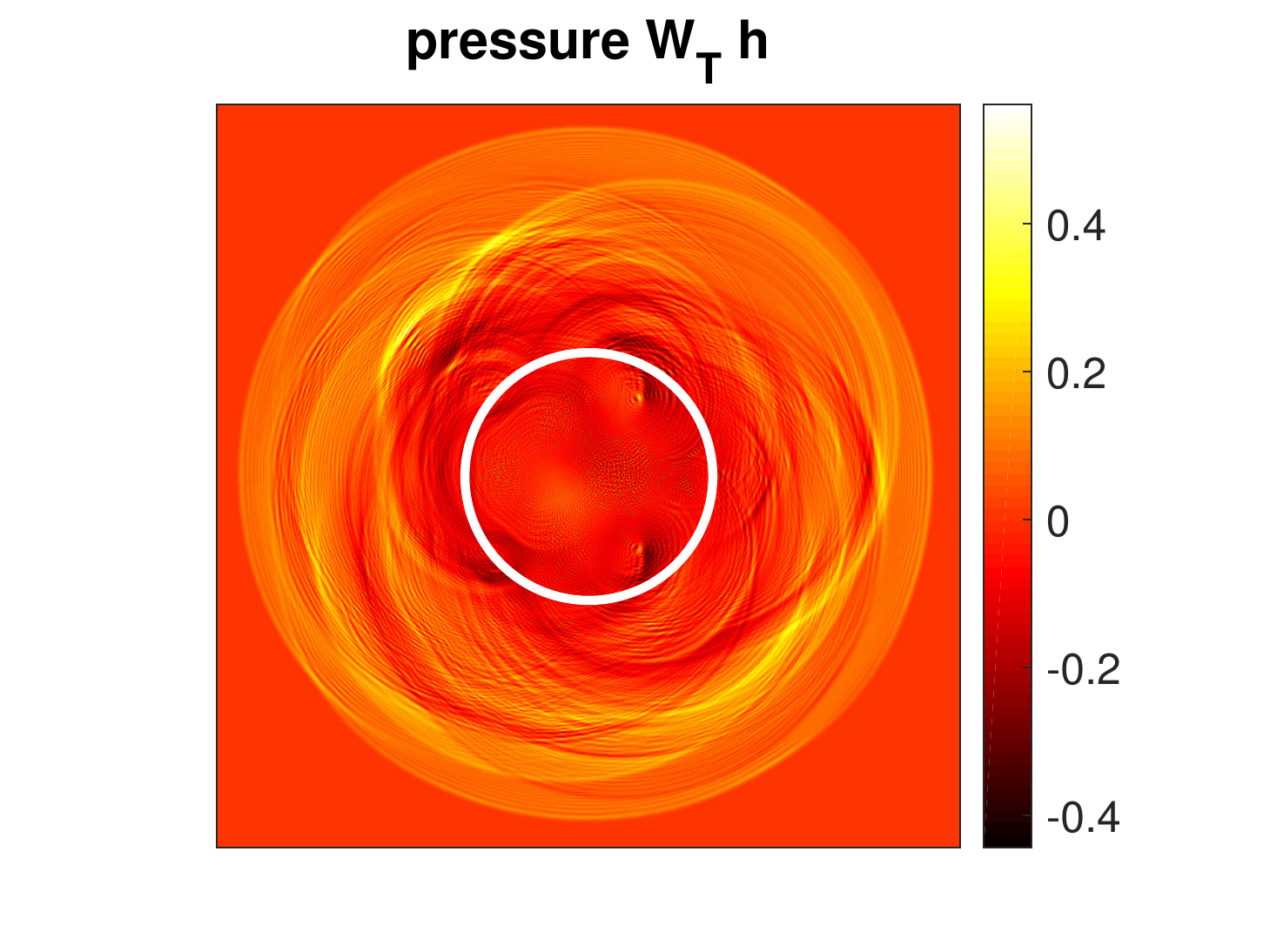} \;
\includegraphics[width=0.48\textwidth]{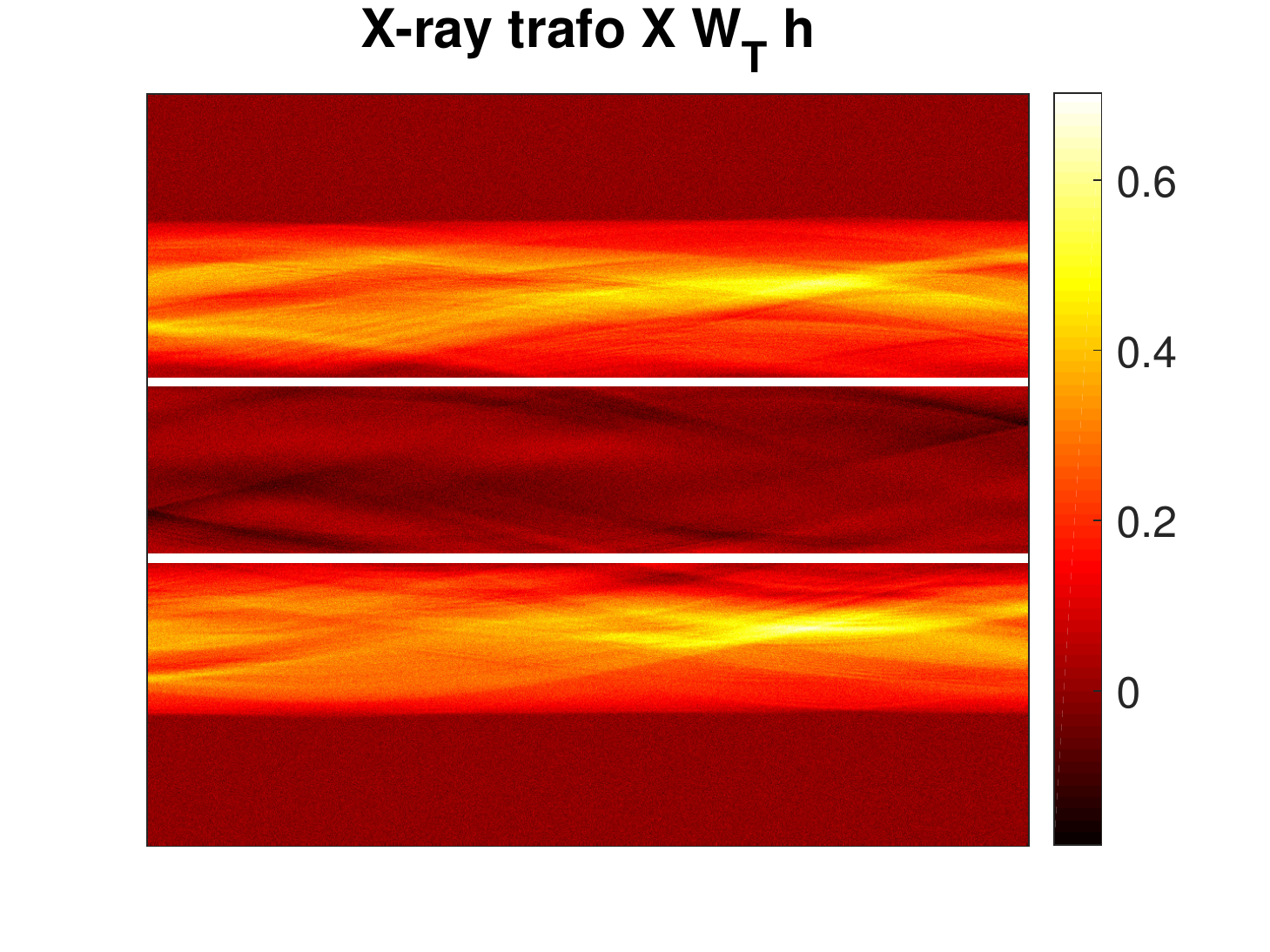} \\[1em]
\includegraphics[width=0.48\textwidth]{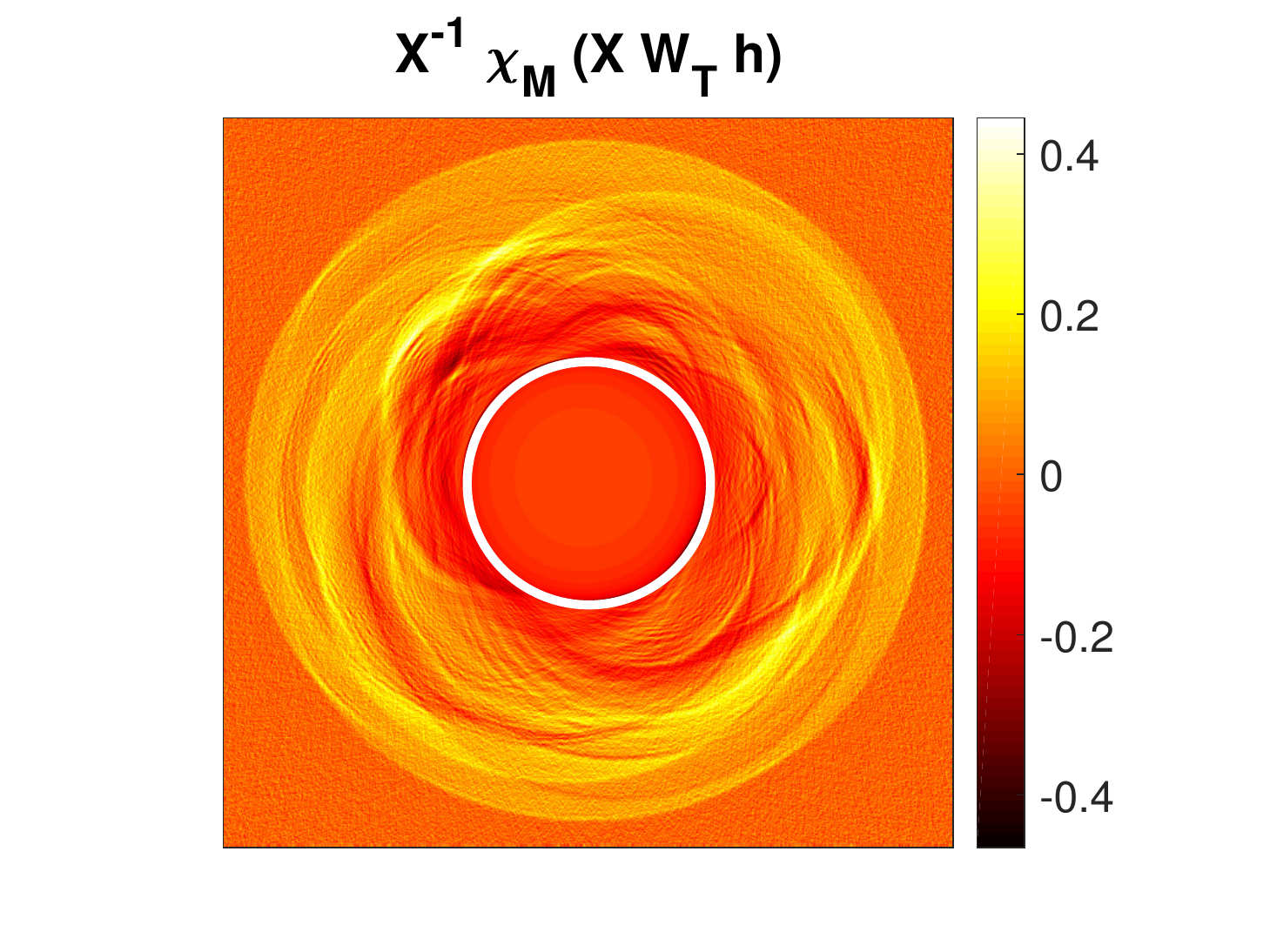}  \;
\includegraphics[width=0.48\textwidth]{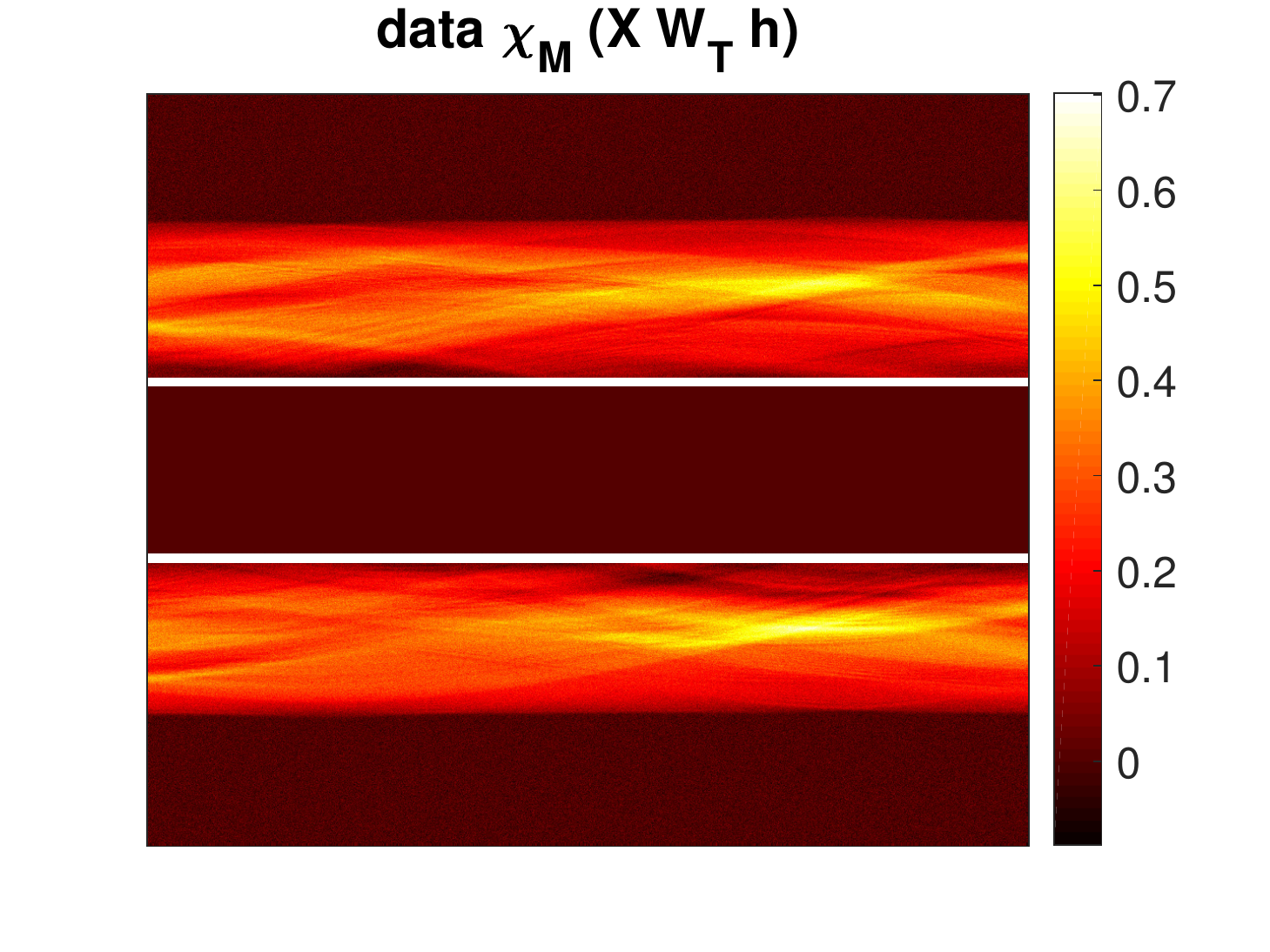}
\caption{ \textbf{Data simulation.}  Top left: Pressure
$\Wo_T \source$ at time $T= 2$.
Top  right:  Noisy X-ray transform
$ \Xo \Wo_T \source + \mathrm{noise}$.
Bottom right: Noisy data $G = \chi_ D  \Xo \Wo_T \source + \mathrm{noise}$ used as input for the one-step methods.
Bottom left:   Recovered final time pressure
$ \Xo^{-1} G  \simeq  \Wo_T \source$ used as input for the second step of the two-step method. The white circles indicate the
imaging domain $D_1$ and  the stripe  inside the white lines
is the missing region $[-1,1]$. } \label{fig:data}
\end{figure}

\subsection{Simulation and reconstruction  results}

Figure~\ref{fig:phantom} (left) shows the PA source to be reconstructed.  We used a trapping sound speed \cite{qian2011new}  shown in Figure~\ref{fig:phantom} (right).
The simulated data are shown in Figure~\ref{fig:data}.
The top left picture  shows the pressure $\Wo_T  \source$
and the top right  picture   shows its X-ray transform $\Xo \Wo_T  \source$  to which  we have added Gaussian white noise with a   standard deviation  of $\SI{20}{\percent}$ of the
mean value of $\Xo \Wo_T \source$, resulting  in a
relative $L^2$-data error $\snorm{\mathrm{noise}}_{L^2} /\snorm{\Xo \Wo_T \source}_{L^2}  \simeq  \SI{11.73}{\percent}$.
 The FFD data $g = \chi_M \Xo \Wo_T  \source + \mathrm{noise}$
 are  shown  in the bottom right image in
 Figure  \ref{fig:data}. The bottom left image shows
 $\Xo^{-1} g$ which is the input for the second step of the two-step method.

\begin{figure}[htb!] \centering
\includegraphics[width=0.48\textwidth]{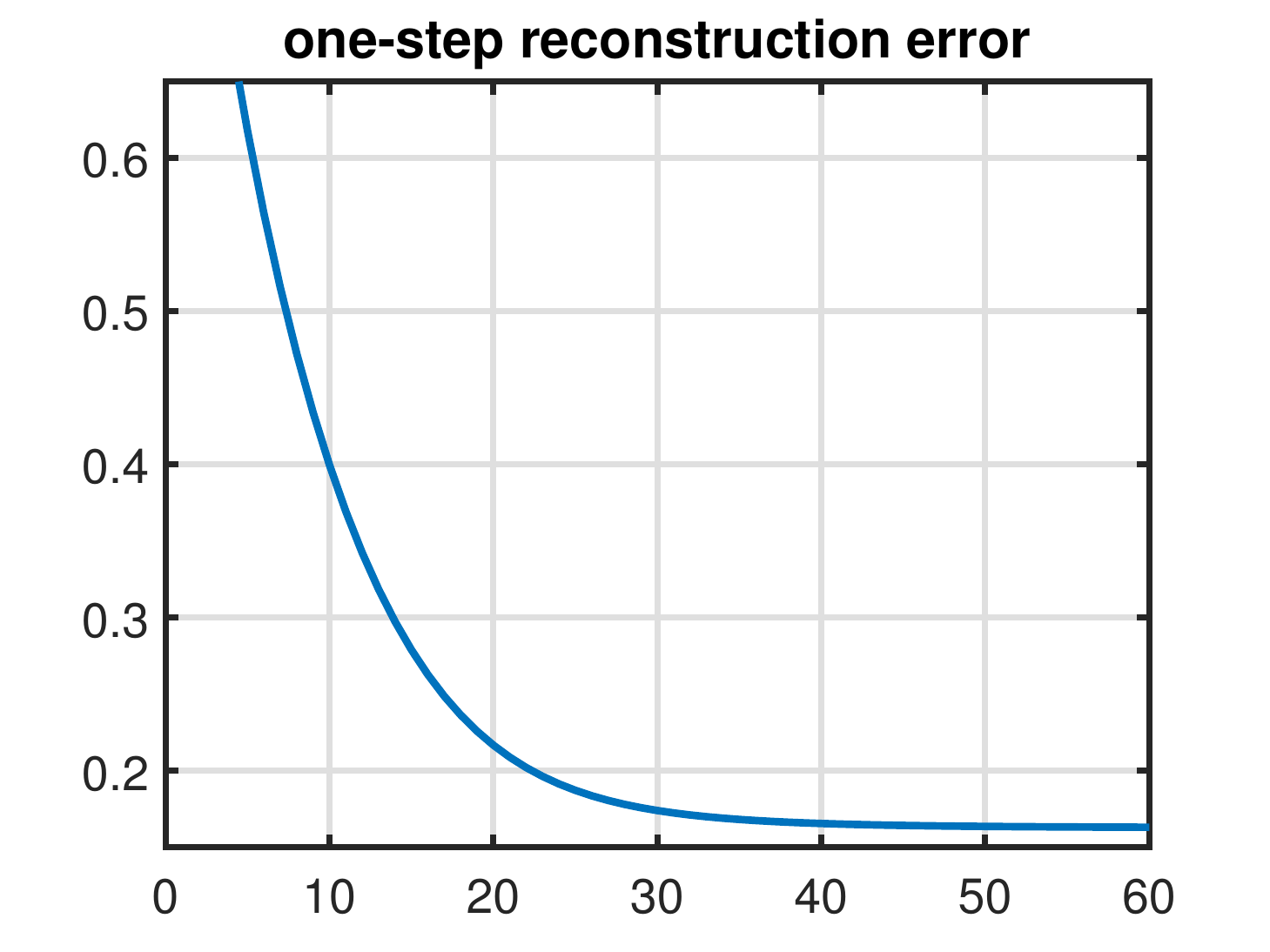} \;
\includegraphics[width=0.48\textwidth]{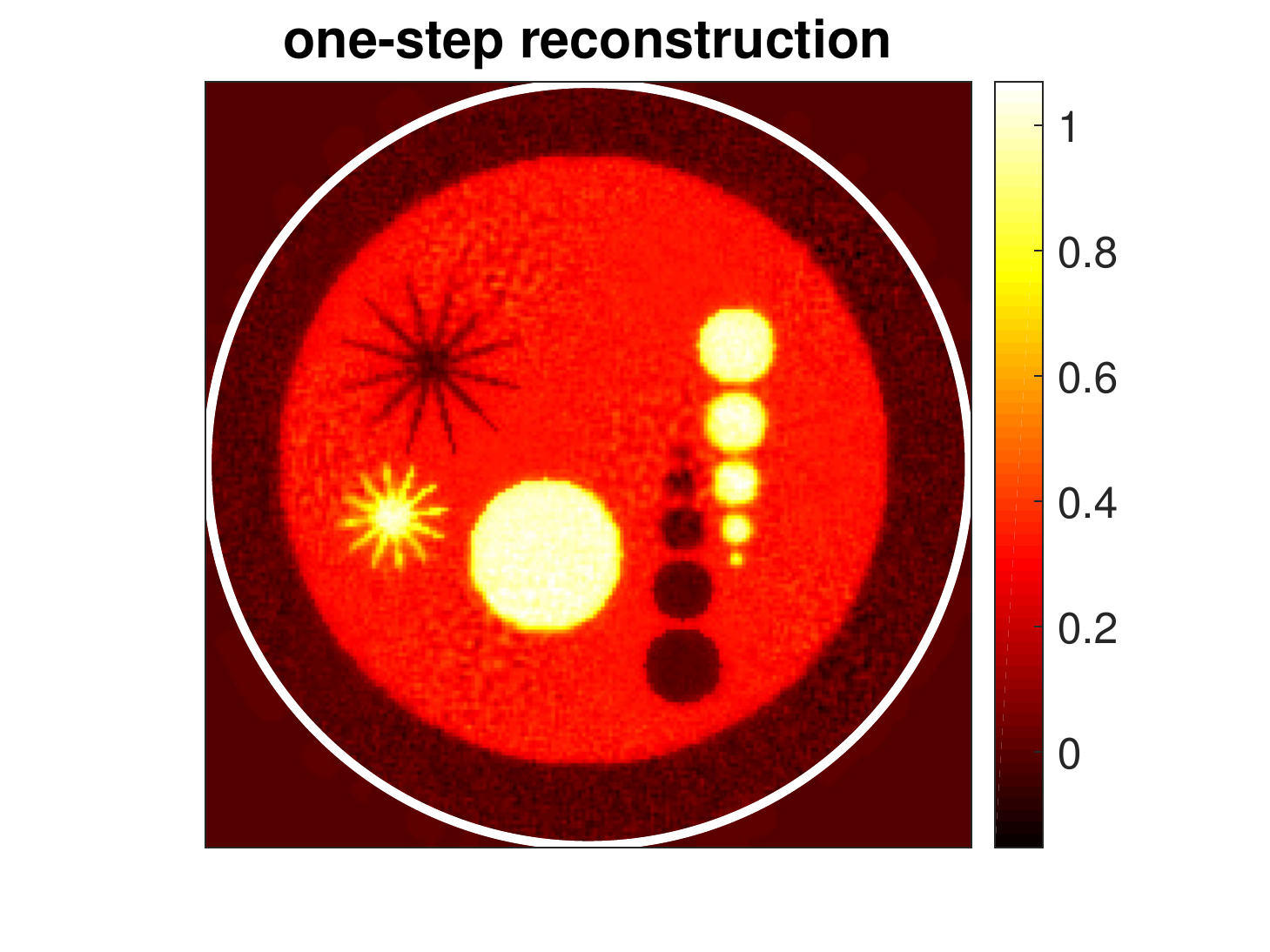} \\[1em]
\includegraphics[width=0.48\textwidth]{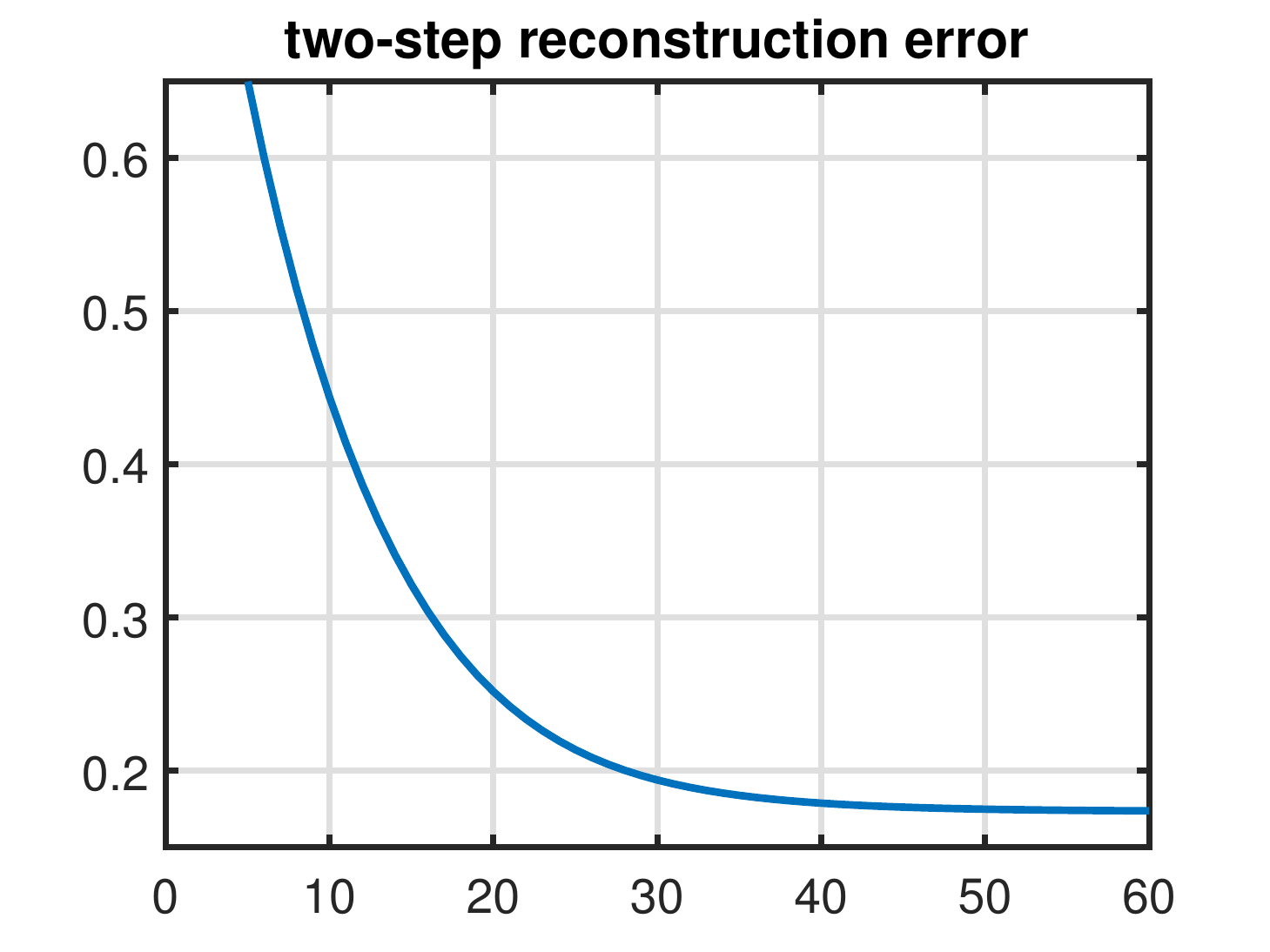} \;
\includegraphics[width=0.48\textwidth]{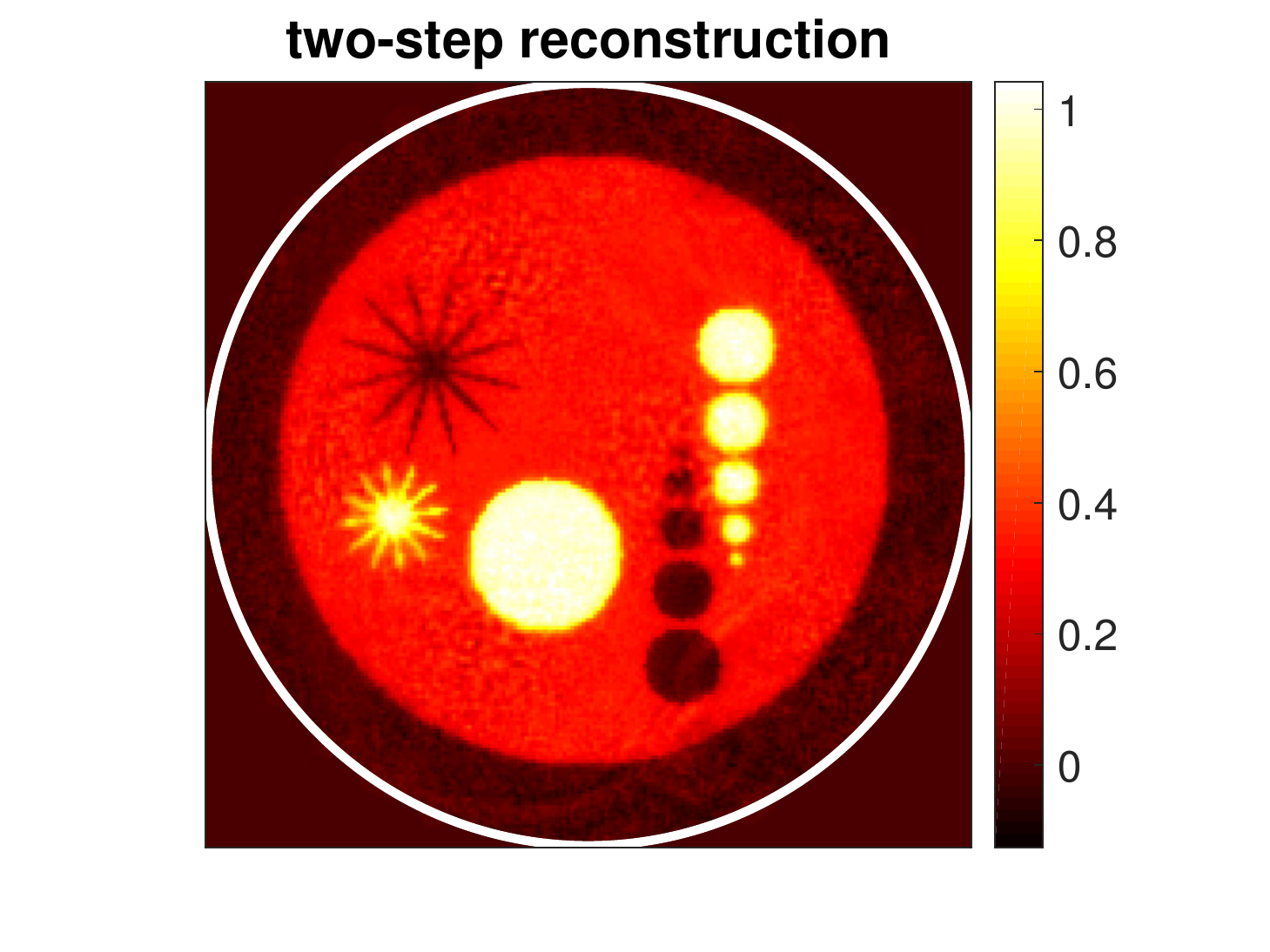}\\[1em]
\includegraphics[width=0.48\textwidth]{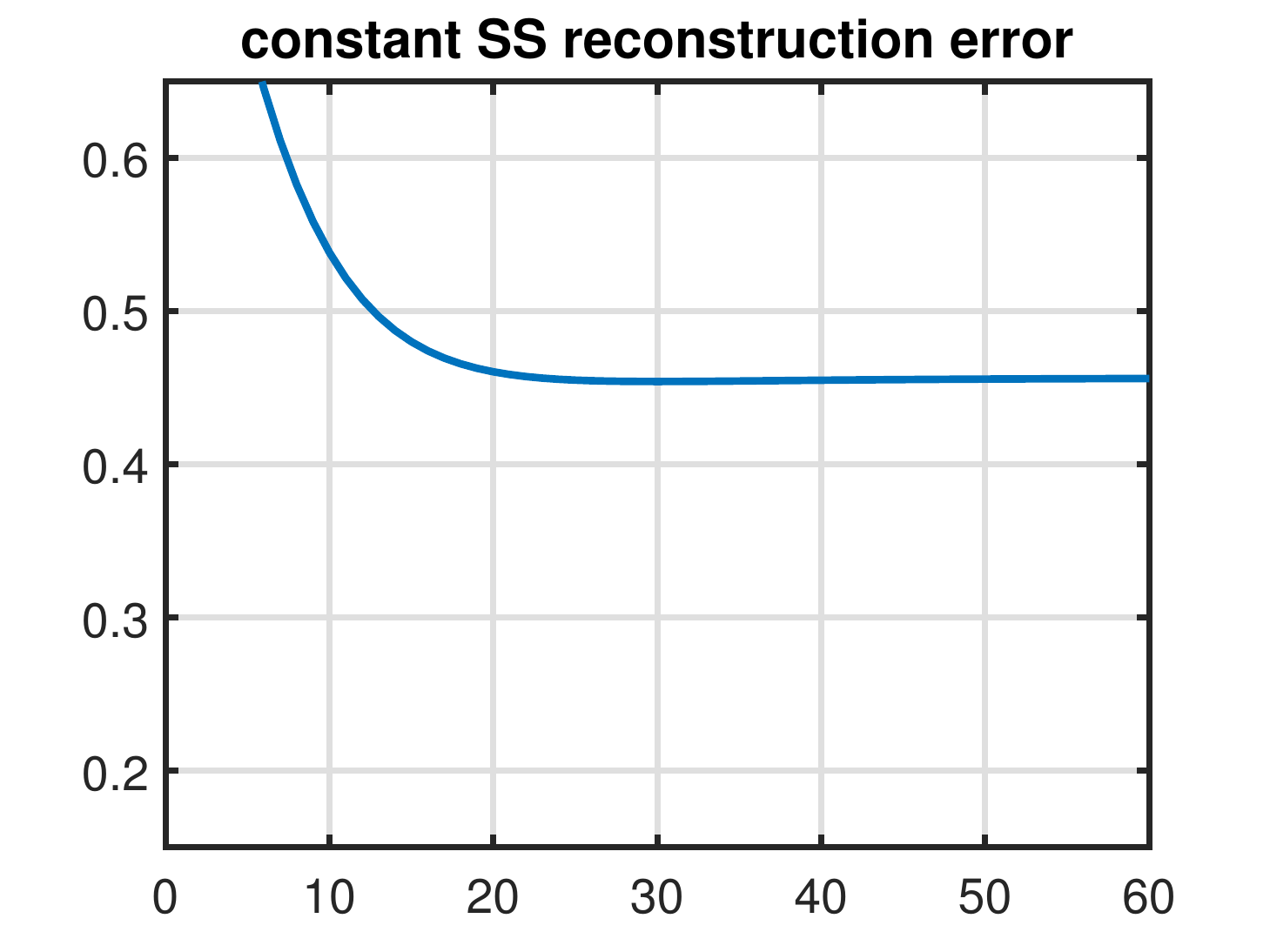} \;
\includegraphics[width=0.48\textwidth]{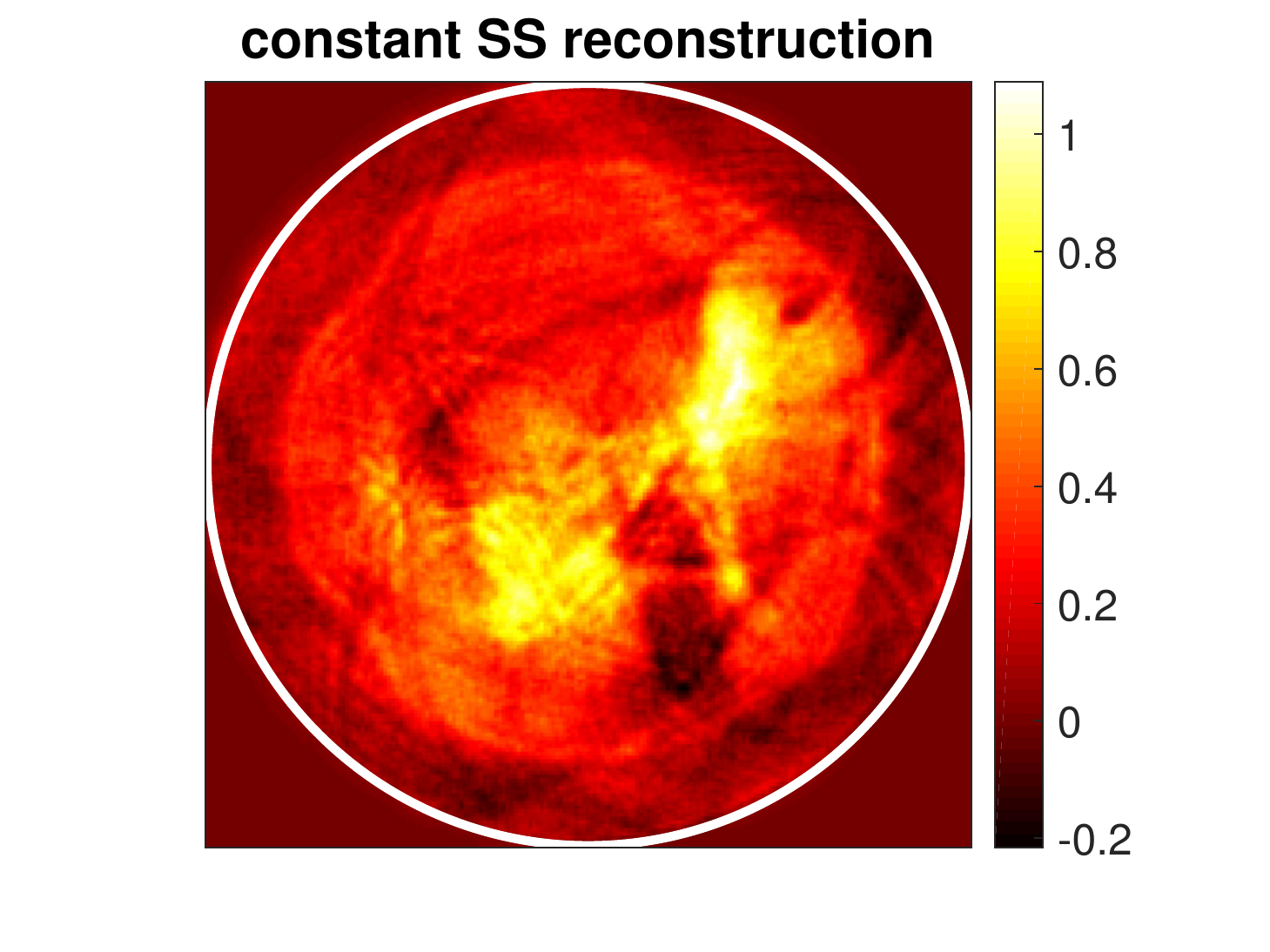}
\caption{ \textbf{Reconstruction results.}
Top left: Relative $L^2$-reconstruction error  with the one-step method.
Top right: Reconstruction with one-step methods after 60 iterations.
Center left: Relative $L^2$-reconstruction error  with the two-step method.
Center right: Reconstruction with two-step methods after 60 iterations.
Bottom left: Relative $L^2$-reconstruction error  assuming 
 constant sound speed for the reconstruction algorithm.
Bottom right: Reconstruction assuming constant sound speed after 60 iterations.}
 \label{fig:recon}
\end{figure}

Reconstruction results with the  one-step (top row) and the two-step (center row) 
method are shown in Figure~\ref{fig:recon}. The images on the left hand side show the  evaluation of  the relative $L^2$-reconstruction error
$\snorm{ \source_k -\source_\star }_{L^2} /\snorm{\source_\star}_{L^2} $  in
 dependence of the iteration index $k$. The images on the
 right show the reconstructions results of the one-step method (top) and two-step method (center) after 60 iterations. The relative $L^2$-reconstruction   error is $\SI{16.33}{\percent}$ for the one-step and
 $\SI{17.39}{\percent}$ for the two-step method. Visually as well as in
 terms of the $L^2$-error, the one-step method shows
 slightly improved results compared to the two-step  method.
 However, as already observed  in  \cite{zangerl2018full}
 the two-step method works remarkably well. The bottom row in Figure~\ref{fig:recon} shows the reconstruction results assuming a constant sound speed  for image reconstruction from the data shown in Figure~\ref{fig:data}.
The constant sound speed has been taken as the average value of the actual sound speed used for data generation. One observes that the reconstruction  results using a constant  sound speed 
are  inferior, demonstrating that it is crucial  to take  
sound speed variations into account.

\section{Conclusion}
\label{sec:conclusion}

In this paper, we have  studied   image reconstruction
using full field detection  PAT. We reviewed existing
result for the  inversion  problem  and introduced   new
preconditioned one-step reconstruction methods with
and without explicit regularization.   We presented numerical results
for the 2D limited  data  setting  using  the preconditioned
forward backward  one-step splitting Algorithm 8.
We compared the results with the previous two-step reconstruction
approach,  where in the first step the final wave data $\Wo_T \source$
are approximately  recovered by  applying the inverse
X-ray transform to the full field data $\chi_M  ( \Xo  \Wo_T  \source)$.
For a fair comparison,  we used  the same
regularization term and minimization  algorithm for the
one-step and the two-step method.
As shown in Figure \ref{fig:recon} the
one-step   as well as the two-step  method produce accurate
results. The one-step method, however,  yields  slightly
fewer visual artefacts   and reduces the relative   $L^2$-reconstruction
error compared to the two-step method.

There are   several open  problems related to the FFD-PAT inversion Problem~1. Uniqueness and stability of reconstruction are  known
for the full data case, but in the case of limited case neither stability nor uniqueness of reconstruction are known. This is also the case
for the related final time  wave inversion Problem~2.
Our numerical investigations suggest  uniqueness as well as stability
in the case $T$ is taken sufficiently  large and the measurement domain is  sufficiently  large. Theoretically investigating these  issues  will be subject of future research. Moreover, numerical and experimental investigations in  3D will be performed, especially for the case where the  measurement domains   in the field data  are small
and where limited  data  artifact are expected.

\section*{Acknowledgments}

M.H.  and G.Z acknowledge support of the Austrian Science Fund (FWF), project P 30747-N32.
The research of L.N. is supported by the National
Science Foundation (NSF) Grants DMS 1212125 and DMS 1616904.
The work of R.N. has been supported by the FWF, project P 28032.

\end{document}